\documentclass{amsart}

\usepackage{amsfonts}
\usepackage{amsmath}
\usepackage{amsthm}
\usepackage{amssymb}
\usepackage[english]{babel}
\usepackage{graphics}
\usepackage{latexsym}
\usepackage{longtable} 
\usepackage{mathrsfs} 
\usepackage{graphicx}
\usepackage{wrapfig} %\usepackage{txfonts}
\usepackage[pdftex,colorlinks=true,linkcolor=blue,citecolor=red]{hyperref}
\usepackage{enumitem}
\usepackage{slashed}
\usepackage{esint}
\usepackage{mathcomp}
\usepackage{stmaryrd}
\usepackage{wasysym}
\usepackage{fancyhdr}

\newtheorem{theorem}{Theorem}
\newtheorem{corollary}[theorem]{Corollary}
\newtheorem{lemma}[theorem]{Lemma}
\newtheorem{proposition}[theorem]{Proposition}

\newtheorem{claim}[theorem]{Claim}
\newtheorem{example}[theorem]{Example}

\theoremstyle{definition}
\newtheorem{definition}[theorem]{Definition}

\renewcommand{\d}{\mathrm{d}}

\newcommand{\mL}{\mathcal{L}}
\newcommand{\mH}{\mathcal{H}}

\newcommand{\D}{\mathrm{D}}

\newcommand{\R}{\mathbb{R}}

\newcommand{\N}{\mathbb{N}}

\newcommand{\mB}{\mathbb{B}}

\newcommand{\ms}{\medskip}

\newcommand{\al}{\alpha}
\newcommand{\be}{\beta}

\newcommand{\Ga}{\Gamma}
\newcommand{\de}{\delta}

\newcommand{\e}{\varepsilon}

\newcommand{\La}{\Lambda}
\newcommand{\ka}{\kappa}
\newcommand{\Om}{\Omega}
\newcommand{\om}{\omega}

\newcommand{\ze}{\zeta}

\newcommand{\av}{-\hspace{-10.5pt}\displaystyle\int}

\newcommand{\weak }{\, -\!\!\!\!\!-\!\!\!\!\rightharpoonup}
\newcommand{\weakstar }{ \overset{\, *_{\phantom{|}}}{{\smash{\, -\!\!\!\!-\!\!\!\!\rightharpoonup}}\, } }

\newcommand{\larrow}{\longrightarrow}
\newcommand{\ot}{\otimes}

\newcommand{\LL}{\text{\LARGE$\llcorner$}}
\newcommand{\p}{\partial}
\newcommand{\sub}{\subseteq}
\newcommand{\set}{\setminus}
\newcommand{\by}{\times}

\renewcommand{\div}{\mathrm{div}}

\newcommand{\bt}{\begin{theorem}}\newcommand{\et}{\end{theorem}}
\newcommand{\bd}{\begin{definition}}\newcommand{\ed}{\end{definition}}
\newcommand{\bl}{\begin{lemma}}\newcommand{\el}{\end{lemma}}
\newcommand{\beq}{\begin{equation}}\newcommand{\eeq}{\end{equation}}
\newcommand{\bc}{\begin{claim}}\newcommand{\ec}{\end{claim}}
\newcommand{\bex}{\begin{example}}\newcommand{\eex}{\end{example}}
\newcommand{\bcor}{\begin{corollary}}\newcommand{\ecor}{\end{corollary}}
\newcommand{\bp}{\begin{proof}}\newcommand{\ep}{\end{proof}}

\newcommand{\BPL}{\medskip \noindent \textbf{Proof of Lemma} }

\newcommand{\BPP}{\medskip \noindent \textbf{Proof of Proposition} }

\numberwithin{equation}{section}
\makeatletter
\@namedef{subjclassname@2020}{\textup{2020} Mathematics Subject Classification}
\makeatother

\begin{document}

\title[Second order vectorial $\infty$-eigenvalue problems]{Generalised second order vectorial $\infty$-eigenvalue problems}

\author{Ed Clark}

\address{E. C., Department of Mathematics and Statistics, University of Reading, Whiteknights Campus, Pepper Lane, Reading RG6 6AX, United Kingdom}

\email{e.d.clark@pgr.reading.ac.uk}
 
\author{Nikos Katzourakis}

\address[Corresponding author]{N. K., Department of Mathematics and Statistics, University of Reading, Whiteknights Campus, Pepper Lane, Reading RG6 6AX, United Kingdom}

\email{n.katzourakis@reading.ac.uk}

\thanks{E.C.\ has been financially supported through the UK EPSRC scholarship GS19-055}
  

\date{}


\begin{abstract} We consider the problem of minimising the $L^\infty$ norm of a function of the hessian over a class of maps, subject to a mass constraint involving the $L^\infty$ norm of a function of the gradient and the map itself. We assume zeroth and first order Dirichlet boundary data, corresponding to the ``hinged" and the ``clamped" cases. By employing the method of $L^p$ approximations, we establish the existence of a special $L^\infty$ minimiser, which solves a divergence PDE system with measure coefficients as parameters. This is a counterpart of the Aronsson-Euler system corresponding to this constrained variational problem. Furthermore, we establish upper and lower bounds for the eigenvalue.
\end{abstract}

\maketitle

%%%%%%%%%%%%%%%%%%%%
\!\!
%%%%%%%%%%%%%%%%%%%

\section{Introduction and main results}   \label{Section1}

Let $n, N \in \N$ with $n\geq 2$, and let $\Om \Subset \R^n$ be a bounded open set with Lipschitz boundary $\partial \Om$. In this paper we are interested in studying nonlinear second order $L^\infty$ eigenvalue problems. Specifically, we investigate the problem of finding a minimising map $u_{\infty}:\overline{\Om}\longrightarrow \R^N$, that solves
\beq
\label{1.1}
\begin{split}
\|f(\D^2u_\infty)\|_{L^{\infty}(\Om)}\,=\, \inf \Big\{& \|f(\D^2v)\|_{L^{\infty}(\Om)}\ :\\
& v\in W^{2,\infty}_{\mathrm B}(\Om;\R^N), \ \|g(v, \D v)\|_{L^{\infty}(\Om)}=1\Big\}.
\end{split}
\eeq 
Additionally, we pursue the necessary conditions that these constrained minimisers must satisfy, in the form of PDEs. In the above, $f: \R^{N \by n^2}_s \longrightarrow \R$ and $g :\R^{N} \times \R^{N\times n} \longrightarrow \R$ are given functions, that will be required to satisfy some natural assumptions, to be discussed later in this section. We note that $\R^{N \by n^2}_s$ symbolises the symmetric subspace of the tensor space $\R^N \ot (\R^n \ot \R^n)$ wherein the hessians of twice differentiable maps $u : \Om \larrow \R^N$ are valued. The functional Sobolev space $W^{2,\infty}_{\mathrm B}(\Om;\R^N)$ appearing above will taken to be either of:
\beq
\label{1.2}
\left\{ \ \ 
\begin{split}
W^{2,\infty}_{\mathrm C}(\Om;\R^N) :&= \, W^{2,\infty}_0(\Om;\R^N),
\\
W^{2,\infty}_{\mathrm H}(\Om;\R^N) :&= \, W^{2,\infty} \cap W^{1,\infty}_0 (\Om;\R^N).\, 
\end{split}
\right. 
\eeq
The space $W^{2,\infty}_{\mathrm C}(\Om;\R^N)$ encompasses the case of so-called clamped boundary conditions, which can be seen as first order Dirichlet or as coupled Dirichlet-Neumann conditions, requiring $|u| =|\D u|= 0$ on $\partial \Om$. On the other hand, $W^{2,\infty}_{\mathrm H}(\Om;\R^N)$ encompasses the so-called hinged boundary conditions, which are zeroth order Dirichlet conditions, requiring $|u| = 0$ on $\partial \Om$. This is standard terminology for such problems, see e.g.\ \cite{KP}.  

Problem \eqref{1.1} lies within the Calculus of Variations in $L^{\infty}$, a modern area, initiated by Gunnar Aronsson in the 1960s. Since then this field has undergone a substantial transformation. There are some general complications one must be wary of when tackling $L^{\infty}$ variational problems. For example, the $L^{\infty}$ norm is generally not Gateaux differentiable, therefore the analogue of the Euler-Lagrange equations cannot be derived directly by considering variations. Any supremal functional also has issues with locality in terms of minimisation on subdomains. Further, the space itself lacks some fundamental functional analytic properties, such as reflexivity and separability. Higher order problems and problems involving constraints present additional difficulties and have been studied even more sparsely, see e.g.\ \cite{AB, BJ, CKM, CK, K1, K2, K3, K4, KM, KPr}. In fact, this paper is an extension of \cite{K4} to the second order case, and generalises part of the results corresponding to the existence of minimisers and the satisfaction of PDEs from \cite{KP}. In turn, the paper \cite{K4} generalised results on the scalar case of eigenvalue problems for the $\infty$-Laplacian (\cite{JL,JLM}). For various interesting results, see for instance \cite{AP, AB, BK, CDP, KZ, MWZ, PZ, RZ}.

The vectorial and higher order nature of the problem we are considering herein precludes the use of standard methods, such as viscosity solutions (see e.g.\ \cite{K0} for a pedagogical introduction). However, we overcome these difficulties by approximating by corresponding $L^p$ problems for finite $p$ case and let $p\to \infty$. The intuition for using this technique is based on the rudimentary idea that, for a fixed $L^\infty$ function on a set of finite measure, its $L^p$ norm tends to its $L^{\infty}$ norm as $p\to \infty$. This technique is rather standard for $L^\infty$ problems, and in the vectorial higher order case we consider herein is essentially the only method known. Even the very recent intrinsic duality method of \cite{BK} is limited to scalar-valued first order problems.

To state our main result, we now introduce the required hypotheses for the functions $f$ and $g$:
\beq
\label{1.3}
 \left\{  
\begin{split}
 & (a) \  f\in C^1(\R^{N \by n^2}_s). \\
& (b) \   f \ \text{is} \ \text{(Morrey) $2$-quasiconvex}.\\
& (c) \ \text{There exist}\ 0< C_1\leq C_2 \text{ such that, }\text{for all}\ X\in \R^{N \by n^2}_s \setminus \{0\},\\
& \ \ \ \ \ \ \ \ \ \ \ \ \ \ \ \ \ \ \ \ \ 0 < C_1f(X)\, \leq \, \partial f(X):X \, \leq \, C_2f(X). \phantom{\Big|}\\
& (d)  \  \text{There exist} \ C_3, ..., C_6>0, \al>1 \text{ and } \be\leq 1: \text{ for all} \ X\in \R^{N \by n^2}_s,\\
&\ \ \ \ \ \ \ \ \ \ \ \ \ \ \ \ \  \ \ -C_3+C_4|X|^{\al}\, \leq \, f(X)\leq C_5|X|^{\al} +C_6, \phantom{\Big|}\\
&\ \ \ \ \ \ \ \ \ \ \ \ \ \ \ \ \  \ \ \ \ \ \ \ \ \ \ \ \   |\partial f(X)|\, \leq \, C_5 f(X)^{\be}+C_6.\phantom{_\big|}
\end{split}
\right. 
\eeq
\beq
\label{1.4}
\left\{ 
\begin{split}
 & (a) \ g\in C^1(\R^{N}\times \R^{N\times n}). \\
& (b) \  g \ \text{is coercive, in the sense that}\\
& \ \ \ \ \ \ \ \ \ \ \ \ \ \ \ \ \ \ \ \ \ \ \ \ \ \ \ \  \lim_{t\to \infty} \bigg(\inf_{(\eta,P)\in \R^N \! \by \R^{N\times n}, |(\eta,P)|=1}g(t\eta, tP)\bigg)=\infty.\phantom{\Big|}
\\
& (c) \ \text{There exist}\ 0<C_7\leq C_8: \text{ for all } (\eta,P)\in \big(\R^N \! \by \R^{N\times n}\big) \! \setminus \! \{(0,0)\},\\
& \ \ \ \ \ \  0<C_7\, g(\eta, P)\,  \leq \,  \p_\eta g (\eta, P)\cdot \eta+\p_P g (\eta, P):P\; \leq \; C_8\,g(\eta, P). \phantom{\Big|}
\end{split}
\right. \ \
\eeq
In the above, $\partial f(X)$ denotes the the derivative of $f$ whilst $\p_\eta g $ and $\p_P g $ signifies the respective partial derivatives. Additionally ``:" and $``\cdot"$ represent the Euclidean inner products. The terminology of (Morrey) $2$-quasiconvex refers to the standard notion for integral functionals for higher order functionals (see e.g.\ \cite{C, DFLM, D, Zh}), namely
\[
\ \ \ F(X) \, \leq \, \av_\Om F(X+\D^2 \phi)\, \mathrm d \mL^n, \ \ \ \ \forall\ \phi \in W^{2,\infty}_0(\Om;\R^N),\  \forall\ X\in \R^{N \by n^2}_s.
\]
We note that herein we will be using the following function space symbolisations:
\[
\begin{split}
C^2_{\mathrm B}(\overline{\Om};\R^N) \, &:=\, C^2(\overline{\Om};\R^N) \cap W^{2,\infty}_{\mathrm B}(\Om; \R^N),
\\
W^{2,p}_{\mathrm C}(\Om;\R^N) &:= \, W^{2,p}_0(\Om;\R^N),\ \ p\in[1,\infty),
\\
W^{2,p}_{\mathrm H}(\Om;\R^N) &:= \, W^{2,p} \cap W^{1,p}_0 (\Om;\R^N), \ \ p\in[1,\infty),
\end{split}
\]
Further, we will be using the rescaled $L^p$ norms for $p\in[1, \infty)$, given by
\[ 
\| h\|_{L^p(\Om)}\,:=\,\bigg(\, \frac{1}{\mL^{n}(\Om)}\int_{\Om}|h|^p \, \d {\mathcal{L}}^{n}\bigg)^{\frac{1}{p}}\,=\, \bigg(\, {\av}_{\!\!\!\Om}|h|^p \, \d {\mathcal{L}}^{n}\bigg)^{\frac{1}{p}}.
\]
Finally, we observe that \eqref{1.3}(c), implies that $f>0$ on $\R^{N \by n^2}_s\setminus\{0\}$, $f(0)=0$ and $f$ is radially increasing, meaning that $t\mapsto f(tX)$ is increasing on $(0, \infty)$ for any fixed $X\in \R^{N \by n^2}_s\setminus\{0\}$. Similarly, \eqref{1.4}(c) implies that $g>0$ on $(\R^N \times \R^{N\times n}) \setminus\{(0,0)\}$, $g(0,0)=0$  and $g$ is radially increasing on $\R^N \times \R^{N\times n}$, namely $t\mapsto g(t\eta,tP)$ is increasing on $(0, \infty)$ for any fixed $(\eta,P)\in (\R^N \times \R^{N\times n})\setminus\{(0,0)\}$.

Below is our main result, in which we consider both cases of boundary conditions simultaneously.

\begin{theorem}\label{1}
Suppose that the  assumptions \eqref{1.3} and \eqref{1.4} hold true. Then: 
\smallskip
\\
(A) The problem \eqref{1.1} has a solution $u_{\infty}\in W^{2,\infty}_{\mathrm B}(\Om;\R^N).$\smallskip\\
(B) There exist Radon measures 
\[
{\mathrm M}_\infty \in \mathcal{M}\big(\overline{\Om}; \R_s^{N \by n^2}\big), \ \ \ \nu_{\infty}\in \mathcal{M}(\overline{\Om}),
\]
such that
\beq
\label{1.5}
\begin{split}
 \int_{\overline{\Om}}   \D^2\phi : \d {\mathrm M}_\infty=\Lambda_{\infty} \int_{\overline{\Om}}
\Big( \p_\eta g (u_{\infty}, \D u_{\infty})\cdot \phi \, +\, \p_P g (u_{\infty}, \D u_{\infty}) : \D \phi \Big) \, \d \nu_{\infty}
\end{split}
\eeq
for all test maps $\phi\in C^2_{\mathrm B}(\overline{\Om}; \R^N)$, where
\beq
\label{1.6}
\begin{split}
 \Lambda_{\infty} \, =\, \big\|f(\D^2(u_{\infty}) \big\|_{L^{\infty}(\Om)}>0.
\end{split}
\eeq
Additionally, we have the following a priori lower bound for the eigenvalue
\beq
\label{1.7}
\Lambda_{\infty}\geq \Bigg ( 
  \frac{  C_4 } {\mathrm{diam}(\Om)^\al \Big(C(\infty,\Om) \|\p_\eta g \|_{L^{\infty}(\{g \leq1\})}+\|\p_P g \|_{L^{\infty}(\{g \leq1\})}\Big)^\al }
 -C_3
 \Bigg)^{\!\!+},
\eeq
where $(\, \cdot\, )^+$ symbolises the positive part, and $C(\infty,\Om)$ equals either the constant of the Poincar\'e inequality (in the case of clamped boundary conditions), or the constant of the Poincar\'e-Wirtinger inequality (in the case of hinged boundary conditions), both taken for $p=\infty$. 
\\
If additionally the boundary $\p\Om$ is $C^2$, we have the a priori upper bound
\beq
\label{1.8A}
\begin{split}
\La_\infty \, \leq & \ C_6 \,+\,  C_5\frac{2^{5\al}}{(c \om(n))^\al} \bigg( \! 1 + \underset{0\leq t \leq 1}{\sup}R(t) \! \bigg)^{\!\!\al} \Big(2^{3n} + \underset{i=1,...,n-1}{\max}\big(\|\ka_i\|_{C^0(\p\Om)}\big)^n  \Big)^{\!\al} \centerdot
\\
&\centerdot \Bigg\{ 1 + \bigg( \! 1+\frac{C}{\e_0^{n+1}} \! \bigg) \mH^{n-1}(\p\Om) \, +\, \sum_{i=1}^{n-1} \left\| \frac{ \ka_i \circ \mathrm{P}_{\Om} }{ 1 - ( \ka_i \circ \mathrm{P}_{\Om})d_\Om } \right\|_{L^\infty(\{d_\Om< \e_0\} \cap \Om)}  \!\!\Bigg\}^{\!\al} ,
\end{split}
\eeq
where $c,C>0$ are dimensionless universal constants, $\om(n)$ is the volume of the unit ball in $\R^n$, $\mH^{n-1}(\p\Om)$ is the perimeter of $\Om$, $\{\ka_1,...,\ka_{n-1}\}$ are the principal curvatures of $\p\Om$, $\mathrm{P}_{\Om}$ is the orthogonal projection on $\p\Om$, $d_\Om$ the distance function of $\p\Om$, $\e_0$ is the largest 
\[
\e \in \left(0 , \min\left\{1\, ,\, \underset{i=1,...,n-1}{\min}\frac{1}{ \|\ka_i\|_{C^0(\p\Om)}} \right\}\right),
\]
for which we have that $d_\Om \in C^2(\{d_\Om\leq \e\} \cap \overline{\Om})$ and $R(t)$ is the smallest radius of the $N$-dimensional ball, for which the sublevel set $\{g\leq t\}$ is contained into the cylinder $\bar\mB^N_{R(t)}(0) \by \R^{N\by n}$, namely
\[
R(t) \,:=\, \inf \Big\{ R>0 \, : \, \{g\leq t\} \sub \mB^N_R(0) \by \R^{N\by n}\Big\}.
\]
\\
(C) The quadruple $(u_{\infty}, \Lambda_{\infty}, {\mathrm M}_\infty,  \nu_{\infty})$ satisfies  the following approximation properties: there exists a sequence of exponents $(p_j)_1^{\infty} \sub (n/ \al)$ where $p_j\to \infty$ as $j\to \infty$, and for any $p$, a quadruple 
\[
(u_p,\Lambda_p, {\mathrm M}_p, \nu_p)\ \in \ W^{2,\alpha p}_{\mathrm B}(\Om;\R^N)\times [0, \infty) \times\mathcal{M}\big(\overline{\Om}; \R_s^{N \by n^2}\big) \times \mathcal{M}(\overline{\Om}),
\] 
such that
\beq
\label{1.8}
\left\{ \ \
\begin{array}{ll}
u_{p} \larrow u_\infty, &  \ \ \text{in } C^1 \big(\overline{\Om};\R^N \big), \smallskip
\\
		\D^2 u_{p} \weak \D^2 u_\infty, &  \ \ \text{in } L^q \big(\Om; \R^{N \by n^2}_s\big), \ \text{for all } q\in (1,\infty), \smallskip
\\
\Lambda_p \larrow \Lambda_{\infty}, &  \ \ \text{in } [0,\infty),  \smallskip
\\
{\mathrm M}_p \weakstar {\mathrm M}_\infty, &  \ \ \text{in }\mathcal{M}\big(\overline{\Om}; \R_s^{N \by n^2}\big),  \smallskip
\\
\nu_p \weakstar \nu_{\infty}, &  \ \ \text{in } \mathcal{M}(\overline{\Om}),  \smallskip
\end{array}
\right.
\eeq
as $p \to \infty$ along $(p_j)_1^{\infty}$. Further, $u_p$ solves the constrained minimisation problem
\beq
\label{1.9}
\|f(\D^2u_p)\|_{L^{p}(\Om)}\,=\, \inf \Big\{ \|f(\D^2v)\|_{L^{p}(\Om)}\ : \  v\in W^{2,\alpha p}_{\mathrm B}(\Om;\R^N), \ \|g(v, \D v)\|_{L^{p}(\Om)}=1\Big\},
\\
\eeq 
and $(u_p, \Lambda_p)$ satisfies
\beq
\label{1.10}
 \left\{ 
\begin{split}
 & {\av}_{\!\!\!\Om} f(\D^2 u_p)^{p-1} \partial f(\D^2u_p): \D^2 \phi \, \d \mL^n 
 \\
 = &\  (\Lambda_p)^p \, {\av}_{\!\!\!\Om} g(u_p, \D u_p)^{p-1}\Big(\p_\eta g (u_p, \D u_p) \cdot \phi + \p_P g (u_p, \D u_p) : \D \phi \Big) \, \, \d {\mathcal{L}}^{n}
\end{split}
\right.
\eeq
for all test maps $\phi\in W^{2,\alpha p}_{\mathrm B}(\Om;\R^N)$. Finally, the measures ${\mathrm M}_p, \nu_p$ are given by
\beq
\label{1.11}
\left\{ \ \ 
\begin{split}
{\mathrm M}_p &= \frac{1}{\mL^n(\Om)}\bigg(\frac{f(\D^2u_p)}{ \Lambda_p}\bigg)^{\! p-1} \partial f(\D^2u_p) \,  \mL^{n} \LL_{\Om},
\\
\nu_p &= \frac{1}{\mL^n(\Om)}\, g(u_p, \D u_p)^{p-1} \, \mL^{n} \LL_{\Om}.
\\
\end{split}
\right. \ \
\eeq
\end{theorem}

We note that one could pursue optimality in Theorem \ref{1}(A) by using $L^\infty$ versions of quasiconvexity, as developed by Barron-Jensen-Wang \cite{BJW2} but adapted to this higher order case, in regards to the existence of $L^\infty$ minimisers. However, for parts (B) and (C) of Theorem \ref{1} regarding the necessary PDE conditions, we do need Morrey $2$-quasiconvexity, as we rely essentially on the existence of solutions to the corresponding Euler-Lagrange equations and the theory of Lagrange multipliers in the finite $p$ case. Further, the measures ${\mathrm M}_\infty, \nu_{\infty}$ \emph{depend} on the minimiser $u_{\infty}$ in a non-linear fashion, hence one more could perhaps symbolise them more concisely as ${\mathrm M}_\infty(u_\infty), \nu_{\infty}(u_\infty)$. Consequently, the significance of these equations is currently understood to be mostly of conceptual value, rather than of computational nature. However, it is possible to obtain further information about the underlying structure of these parametric measure coefficients. This requires techniques such as measure function pairs and mollifications up to the boundary as in \cite{CK, H, K4}, but to keep the presentation as simple as possible, we refrain from pursuing this -considerably more technical- endeavour, which also requires stronger assumptions.

\section{Proofs}   \label{Section2}

In this section we establish Theorem \ref{1}. Its proof is not labeled explicitly, but will be completed by proving a combination of smaller subsidiary results, including a selection of lemmas and propositions.

Before introducing the approximating problem (the $L^p$ case for finite $p$), we need to establish a convergence result, which shows that the admissible classes of the $p$-problems are non-empty. It is required because the function $g$ appearing in the constraint is not assumed to be homogeneous, therefore a standard scaling argument does not suffice.

\begin{lemma}\label{lemma2}
For any $v\in W^{2,\infty}_{\mathrm B}(\Om;\R^N)\setminus \{0\}$, there exists $(t_p)_{p\in (n/\al, \infty]}$ with $t_p\to t_{\infty}$ as $p\to\infty$, such that
\[ 
\big\|g\big(t_pv, t_p\D v\big)\big\|_{L^p(\Om)}=1, \phantom{\Big|}
\]
for all $p\in (n/\al, \infty]$. Further, if $\|g(v, \D v) \|_{L^{\infty}(\Om)}=1$, then $t_{\infty}=1$.
\end{lemma}

\BPL \ref{lemma2}.
Fix $v\in W^{2,\infty}_{\mathrm B}(\Om;\R^N)\setminus \{0\}$ and define
\[
\rho_{\infty}(t):=\max_{x\in \overline{\Om}}g\big(tv(x),t\D v(x)\big), \ \ \ \ t\geq 0.
\]
It follows that $\rho_{\infty}(0)=0$ and $\rho_{\infty}$ is continuous on $[0, \infty)$. We will now show that $\rho_{\infty}$ is strictly increasing. We first show it is non-decreasing. For any $s>0$ and $(\eta,P) \in \R^N \by \R^{N\times n} \setminus \{(0,0)\}$, our assumption \eqref{1.4}(c) implies
\[
\begin{split}
0\, &<\,\frac{C_7g(s\eta, sP)}{s} 
\\
&\leq \, \p_\eta g (s\eta, sP)\cdot \eta \, +\, \p_P g (s\eta, sP):P
\\
&= \, \p_{(\eta,P)} g (s\eta, sP) : (\eta,P)
\\
&= \, \frac{\d}{\d s}\big(g(s\eta, sP)\big),
\end{split}
\]
thus $s\mapsto g(s\eta, sP)$ is increasing on $(0, \infty)$. Hence, for any $x\in \overline{\Om}$ and $t>s\geq0$ we have $g(tv(x), t \D v(x)) \geq g(sv(x), s\D v(x))$, which yields,
\[
\begin{split}
\rho_{\infty}(s)=\max_{x\in \overline{\Om}}g\big(sv(x),s\D v(x)\big)\, \leq \, \max_{x\in \overline{\Om}}g\big(tv(x),t\D v(x)\big)=\rho_{\infty}(t).
\end{split}
\]
We proceed to demonstrate that $t\mapsto \rho_{\infty}(t)$ is actually strictly monotonic over $(0,\infty)$. Fix $t_0>0$. By Danskin's theorem \cite{Dan}, the derivative from the right $\rho'(t_0^+)$ exists, and is given by the formula
\[
\begin{split}
\rho_{\infty}'(t_0^+) \, = \, \max_{x\in {\Om}_{t_0}}\Big\{ \p_{(\eta,P)} g (t_0v(x), t_0 \D v(x)) : \big(v(x),\D v(x)\big)  \Big\},
\end{split}
\]
where
\[
\begin{split}
\Om_{t_0}\,:=\, \Big \{ \overline{x}\in \overline{\Om}\ :\ \rho_{\infty}(t_0)=g \big(t_0v(\overline{x}), t_0 \D v(\overline{x}) \big)\Big \}.
\end{split}
\]
Hence, by \eqref{1.4}(c) we estimate
\[
\begin{split}
\rho_{\infty}'(t_0^+)&=\frac{1}{t_0}\max_{x\in {\Om}_{t_0}}\Big\{ \p_{(\eta,P)} g (t_0v(x), t_0 \D v(x)) : \big(t_0v(x),t_0\D v(x)\big)  \Big\}
\\
&\geq \frac{C_7}{t_0}\max_{x\in {\Om}_{t_0}}g \big(t_0v(x), t_0\D v(x)\big)\\
&= \frac{C_7}{t_0} \rho_{\infty}(t_0)\\
&>0.
\end{split}
\]
This implies that $\rho_{\infty}$ is strictly increasing on $(0,\infty)$. Next, recall that $g$ is coercive by assumption \eqref{1.4}(b), namely $g(s\eta, sP)\to \infty$ as $s\to \infty$, for fixed $(\eta,P)\neq (0,0)$. Thus, for any fixed point $\overline{x}\in \Om$ with $(v(\overline{x}),\D v(\overline{x}))\ne (0,0)$, which exists because by assumption $v\not \equiv 0$, we have
\[
\begin{split}
\lim_{t\to \infty}\rho_{\infty}(t)\geq\lim_{t \to \infty} g(tv(\overline{x}), t \D v(\overline{x}))=\infty.
\end{split}
\]
Since $\rho_\infty(0)=0$ and $\rho_{\infty}(t) \to\infty$ as $t\to\infty$, by continuity and the intermediate value theorem, there exists a number $t_\infty>0$ such that $\rho_{\infty}(t_{\infty})=1$, that is
\[
\begin{split}
\big\|g\big(t_{\infty}v, t_{\infty} \D v\big)\big\|_{L^{\infty}(\Om)}=1.
\end{split}
\]
If $\|g(v, \D v)\|_{L^{\infty}(\Om)}=1$, then $t_{\infty}=1$, as a result of the strict monotonicity of $\rho_{\infty}$. Now let us fix $p\in (n/\al, \infty)$ and define the continuous function
\[
\begin{split}
\rho_p(t) \, :=\, {\av}_{\!\!\!\Om} g(tv, t \D v)^p \, \d \mL^n , \ \ \ t\geq 0.
\end{split}
\]
Since $g(0,0)=0$, it follows that $\rho_p(0)=0$ and that 
\[
\rho_p(t) \, =\frac{1}{\mL^n(\Om)}\int_{\{(v,\D v)\ne (0,0)\}} g(tv, t\D v)^p \, \d \mL^n.
\]
By Morrey's theorem and our assumptions, we have that $v\in C^1(\overline{\Om};\R^N)\set\{0\}$, therefore $\mL^n\big( \{(v,\D v)\ne (0,0)\} \big)>0$. Consider the family of functions $\{g(tv, t\D v)^p \}_{t>0}$, defined on $\{(v,\D v)\ne (0,0)\}\sub \Om$. By the monotonicity of $s\mapsto g(s\eta, sP)$ on $(0,\infty)$ for $(\eta,P)\neq (0,0)$, for $s<t$ we have 
\[
\text{$g(sv, s\D v)^p \leq g(tv, t\D v)^p$, \ on $\{(v,\D v)\ne (0,0)\}$}. 
\]
Since $g(tv, t\D v)^p \to\infty$ pointwise on $\{(v,\D v)\ne (0,0)\}$ as $t \to\infty$, by the monotone convergence theorem, we infer that 
\[
\int_{\{(v,\D v)\ne (0,0)\}} g(tv, t\D v)^p \, \d \mL^n \larrow \infty,
\]
as $t\to\infty$. As a consequence, $\rho_p(t)\to \infty$ as $t\to \infty$. Since $\rho_p(0)=0$, by the intermediate value theorem there exists $t_p>0$ such that $\rho_p(t_p)=1$, namely 
\[
\big\|g(t_pv, t_p\D v)\big\|_{L^p(\Om)}=1.
\]
For the sake of contradiction, suppose that $ t_p \not\to  t_{\infty}$, as $p\to \infty$. In this case, there exists a subsequence $(t_{p_j})_1^{\infty}\sub (n/\al, \infty)$ and $t_0\in [0, t_{\infty}) \cup (t_{\infty}, \infty]$ such that $t_{p_j}\to t_0$ as $j\to \infty$. Further, $(t_{p_j})_1^{\infty}$ can assumed to be either monotonically increasing or decreasing. We first prove that $t_0$ is finite. If $t_0=\infty$, then the sequence $(t_{p_j})_1^{\infty}$ can be selected to be monotonically increasing. Therefore, by arguing as before, $g(t_{p_j}v, t_{p_j} \D v)\nearrow \infty$ as $j\to \infty$, pointwise on $\{(v,\D v)\ne (0,0)\}$, and the monotone convergence theorem provides the contradiction
\[
1 \, =\,  \lim_{j\to \infty}{\av}_{\!\!\!\Om} g(t_{p_j}v, t_{p_j} \D v)^{p_j} \ \d \mL^n\, =\, {\av}_{\!\!\!\Om} \lim_{j\to \infty} g(t_{p_j}v, t_{p_j} \D v)^{p_j} \ \d \mL^n\, =\, \infty.
\]
Consequently, we have that $t_0\in[0, t_{\infty}) \cup (t_{\infty}, \infty)$. Since $(t_{p_j}v, t_{p_j}\D v)\to (t_0 v, t_0 \D v)$ uniformly on $\overline{\Om}$ as $j\to \infty$, we calculate
\[
\begin{split}
1 \, & =\, \big\|g(t_{p_j}v, t_{p_j}\D v)\big\|_{L^{p_j}(\Om)}
\\
&=\,\big\|g(t_0v, t_0\D v)\big\|_{L^{p_j}(\Om)}\, +\, \mathrm o(1)_{j\to\infty} 
\\
&=\, \big\|g(t_0v, t_0\D v)\big\|_{L^{\infty}(\Om)} \, +\,\mathrm o(1)_{j\to\infty}
\\
&=\, \rho_{\infty}(t_0) \, +\,\mathrm o(1)_{j\to\infty}.
\end{split}
\]
By passing to the limit as $j\to \infty$, we obtain a contradiction if $t_\infty\ne t_0$, because $\rho_{\infty}$ is a strictly increasing function and $\rho_{\infty}(t_\infty)=1$. In conclusion, $t_p\to t_{\infty}$ as $p\to\infty$.
\qed
\ms

Utilising the above result we can now show existence for the approximating minimisation problem for $p<\infty$.

\begin{lemma}\label{lemma3}
For any $p>n/ \al$, the minimisation problem \eqref{1.9} has a solution $u_p \in W^{2,\alpha p}_{\mathrm B}(\Om;\R^N)$.
\end{lemma}

\BPL \ref{lemma3}. Let us fix $p\in(n/ \al, \infty)$ and $v_0 \in W^{2,\infty}_{\mathrm B}(\Om;\R^N)$ where $v_0\, \slashed{\equiv}\, 0$. By application of  Lemma \ref{lemma2}, there exists $t_p>0$ such that $\|g(t_p v_0, t_p\D v_0)\|_{L^p(\Om)}=1$ implying that $t_p v_0$ is indeed an element of the admissible class of  \eqref{1.9}. Hence, we deduce that the admissible class is non empty. Further, by assumption \eqref{1.3}(b), $f$ is (Morrey) $2$-quasiconvex. We now confirm that $f^p$ is also (Morrey) $2$-quasiconvex function, as a consequence of Jensen's inequality: for any fixed $X\in \smash{\R^{N \by n^2}_s}$ and any $\phi \in W^{2,\infty}_0(\Om;\R^N)$, we have
\[ 
f^p(X)\leq \bigg( \, {\av}_{\!\!\!\Om}f(X+\D^2 \phi) \ \d \mL^n \!\bigg)^{\!p}\leq \, {\av}_{\!\!\!\Om}f(X+\D^2 \phi)^p \ \d \mL^n. 
\]
By assumption by assumption \eqref{1.3}(d), we have for some new $C_5(p),C_6(p)>0$ that
\[
f(X)^p \, \leq \, C_5(p)|X|^{\al p}+C_6(p),
\]
for any $X\in \smash{\R^{N \by n^2}_s}$. Moreover, by \cite[Theorem 3.6]{Zh} we have that the functional $v \mapsto \|f(\D^2v)\|_{L^{p}(\Om)}$ is weakly lower semi-continuous on $W^{2,\al p}(\Om;\R^N)$ and therefore the same is true over the closed subspace $W^{2,\al p}_{\mathrm B}(\Om;\R^N)$. Let $(u_i)_{1}^{\infty}$ be a minimising sequence for \eqref{1.9}. As $f\geq 0$, it is clear that $\inf_{i\in \N} \|f(\D^2u_i)\|_{L^{p}(\Om)}\geq 0$. Since the admissible class is non-empty, the infimum is finite. Additionally, by \eqref{1.3}(d), we have the bound
\[
\begin{split}
\inf_{i\in \N} \|f(\D^2u_i)\|_{L^{p}(\Om)} &\leq \big\|f\big(\D^2(t_pv_0)\big)\big\|_{L^{p}(\Om)}
\\
&\leq \big\|C_5\big|t_p\D^2 v_0\big|^{\al}+C_6\big\|_{L^{\infty}(\Om)}
\\
&\leq C_5(t_p)^\al \|\D^2v_0\|^\al_{L^{\infty}(\Om)}+C_6
\\
&< \infty.
\end{split}
\]
Now we show that the functional is coercive in $W^{2,\al p}_{\mathrm B}(\Om;\R^N)$, arguing separately for either case of boundary conditions. By assumption \eqref{1.3}(d) and the Poincar\'e inequality, for any $u\in W^{2,\al p}_{\mathrm C}(\Om;\R^N)$ (satisfying $|u|=|\D u|=0$ on $\p\Om$), we have 
\[
\begin{split}
\bigg( \, {\av}_{\!\!\!\Om} \big|f(\D^2u)+C_3 \big|^p \, \d \mL^n\bigg)^{\!\frac{1}{p}}
\geq C_4\bigg( \, {\av}_{\!\!\!\Om}|\D^2 u|^{\al p} \, \d \mL^n \bigg) ^{\!\frac{1}{p}}
\ \geq C_4'\|u\|^{\al}_{W^{1,\al p}(\Om)},
\end{split}
\]
for a new constant $C_4'=C_4(p)>0$. Hence, for any $u\in W^{2,\al p}_{\mathrm C}(\Om;\R^N)$,
\beq
\label{coercivity}
\|f(\D^2 u)\|_{L^p(\Om)} \, \geq \, C_4'\big(\|u\|_{W^{2,\al p}(\Om)}\big)^{\al} -C_3.
\eeq
The above estimate is also true when $u\in W^{2,\al p}_{\mathrm H}(\Om;\R^N)$, but since in this case we have only $|u|=0$ on $\p\Om$, it requires an additional justification. By the Poincar\'e-Wirtinger inequality involving averages, for any $u\in W^{2,\al p}_{\mathrm H}(\Om;\R^N)$ we have 
\[ 
\left\|\D u- {\av}_{\!\!\!\Om} \D u \ \d \mL^n \right\|_{L^{\al p}(\Om)} \,\leq \, C\|\D^2u\|_{L^{\al p}(\Om)}, 
\]
where $C=C(\al,p,\Om)>0$ is a constant. Since $|u|=0$ on $\p\Om$, by the Gauss-Green theorem we have
\[
\int_{\Om} \D u  \ \d \mL^n= \int_{\partial \Om} u\otimes \hat{n} \ \d \mathcal{H}^{n-1}=0,
\]
where $\mathcal{H}^{n-1}$ denotes the $(n-1)$-dimensional Hausdorff measure. In conclusion, 
\[ 
\big\|\D u \|_{L^{\al p}(\Om)} \, \leq \, C \|\D^2u\|_{L^{\al p}(\Om)}, 
\]
for any $u\in W^{2,\al p}_{\mathrm H}(\Om;\R^N)$. The above estimate together with the standard Poincar\'e inequality applied to $u$ itself allow to infer that \eqref{coercivity} holds for any $u\in W^{2,\al p}_{\mathrm B}(\Om;\R^N)$ in both cases of boundary conditions. Returning to our minimising sequence, by standard compactness results, exists $u_p \in  W^{2,\al p}_{\mathrm H}(\Om;\R^N) $ such that $ u_i \weak u_p$ in $W^{2,\al p}_{\mathrm B}(\Om;\R^N)$, as $i\to\infty$ along a subsequence of indices. Additionally, by the Morrey estimate we have that $u_i\larrow u_p$ in $ C^1(\overline{\Om}; \R^N)$ as $i\to\infty$, along perhaps a further subsequence. Since $u\mapsto \|g(u, \D u)\|_{L^p(\Om)}$ is weakly continuous on $W^{2,\al p}_{\mathrm B}(\Om;\R^N)$, the admissible class is weakly closed in $W^{2,\al p}(\Om;\R^N)$ and hence we may pass to the limit in the constraint. By weak lower semicontinuity of the functional, it follows that a minimiser $u_p$ which satisfies \eqref{1.9} does indeed exist.
\qed
\ms

Now we describe the necessary conditions (Euler-Lagrange equations) that approximating minimiser $u_p$ must satisfy. These equations will involve a Lagrange multiplier, emerging from the constraint $\|g(\cdot, \D (\cdot))\|_{L^p(\Om)}=1$.

\begin{lemma}\label{lemma4}
For any $p>n/ \al$, let $u_p$ be the minimiser of \eqref{1.9} procured by Lemma  \ref{lemma3}.  Then, there exists $\lambda_p \in \mathbb{R}$ such that the pair $(u_p, \lambda_p)$ satisfies the following PDE system
\[
\begin{split}
 &\, {\av}_{\!\!\!\Om} f(\D^2 u_p)^{p-1} \partial f(\D^2u_p): \D^2 \phi \, \d \mL^n 
 \\
&= \,  \lambda_p \, {\av}_{\!\!\!\Om} g(u_p, \D u_p)^{p-1}\Big(\p_\eta g (u_p, \D u_p) \cdot \phi \, +\, \p_P g (u_p, \D u_p) : \D \phi \Big) \, \d {\mathcal{L}}^{n},
\end{split}
\]
for all test maps $\phi\in W^{2,\alpha p}_{\mathrm B}(\Om;\R^N).$
\end{lemma}

In particular, it follows that in both cases $u_p$ is a weak solution in $W^{2,\alpha p}(\Om;\R^N)$ to
\beq
\label{PDEinLp}
\left\{ \ \ 
\begin{split}
 &\, \D^2 : \Big(f(\D^2 u_p)^{p-1} \partial f(\D^2u_p)\Big)
 \\
&= \,  \lambda_p \, \bigg[ g(u_p, \D u_p)^{p-1}\p_\eta g (u_p, \D u_p)  \, -\, \div \Big(g(u_p, \D u_p)^{p-1}\p_P g (u_p, \D u_p) \Big) \bigg],
\end{split}
\right.
\eeq
where we have used the notation $\D^2 : F = \sum_{i,j=1}^n \D^2_{ij} F_{ij}$, when $F \in C^2(\Om;\R^{n \by n})$, which is equivalent to the double divergence (applied once column-wise and once row-wise). Note that in the case of hinged boundary data, we have an additional natural boundary condition arising (since $\D u$ is free on $\p\Om$), we we will not make an particular use of this extra information in the sequel, therefore we refrain from discussing it explicitly.

\BPL \ref{lemma4}. The result follows by standard results on Lagrange multipliers in Banach spaces (see e.g.\ \cite[p.\ 278]{Z}), by utilising assumption \eqref{1.3}(d), which guarantees that the functional is Gateaux differentiable.
\qed
\ms

Now we establish some further results regarding the family of eigenvalues.

\begin{lemma}\label{lemma5} Consider the family of pairs of eigenvectors-eigenvalues $\{(u_p,\lambda_p)\}_{p>n/ \al}$, given by Lemma \ref{lemma4}. Then, for any $p>n/\al$, there exists $\Lambda_p>0$ such that
\[
\lambda_p=\big(\Lambda_p\big)^p >0.
\]
Further, by setting
\[
L_p \, :=\, \big\|f(\D^2 u_p)\big\|_{L^p(\Om)},
\]
we have the bounds
\[
0<\, \bigg(\frac{C_1}{C_8}\bigg)^{\frac{1}{p}}L_p \,\leq \, \Lambda_p \, \leq \, \bigg(\frac{C_2}{C_7}\bigg)^{\frac{1}{p}}L_p.
\]
\end{lemma}

\BPL \ref{lemma5}. We begin by showing that $L_p>0$, namely the infimum over the admissible class of the $p$-approximating minimisation problem is strictly positive, owing to the constraint and our assumptions \eqref{1.3}-\eqref{1.4}. Indeed, there is only one map $u\in W^{2,\al p}_{\mathrm B}(\Om;\R^N)$ for which $\|f(\D^2u)\|_{L^p(\Om)}=0$, namely $u_0\equiv 0$, but this is not an element of the admissible class since $\|g(u_0, \D u_0)\|_{L^p(\Om)}=0$. Now consider the Euler-Lagrange equations in Lemma \ref{lemma4} and select $\phi:=u_p$, to obtain
\[
\begin{split}
 & \,{\av}_{\!\!\!\Om} f(\D^2 u_p)^{p-1} \partial f(\D^2u_p): \D^2 u_p \, \d \mL^n 
 \\
&= \,  \lambda_p \, {\av}_{\!\!\!\Om} g(u_p, \D u_p)^{p-1}\Big(\p_\eta g (u_p, \D u_p) \cdot u_p + \p_P g (u_p, \D u_p) : \D u_p \Big) \, \d {\mathcal{L}}^{n}.
\end{split}
\]
As $f,g\geq 0$ we can manipulate the respective assumptions \eqref{1.3}(c) and \eqref{1.4}(c) to produce the following bounds:
\[
\begin{split}
 C_1{\av}_{\!\!\!\Om} f(\D^2 u_p)^{p}\, \d \mL^n  \,& \leq  \, {\av}_{\!\!\!\Om}f(\D^2u_p)^{p-1} \partial f(\D^2u_p): \D^2 u_p \, \d \mL^n  
 \\
 &\leq  \,  C_2{\av}_{\!\!\!\Om} f(\D^2 u_p)^{p}\, \d \mL^n,
 \end{split}
\]
\[
\begin{split}
 C_7 \, {\av}_{\!\!\!\Om} g(u_p, \D u_p)^{p}\, \d \mL^n
 &\leq \, {\av}_{\!\!\!\Om} g(u_p, \D u_p)^{p-1} \Big(\p_\eta g (u_p, \D u_p) \cdot u_p \,+
 \\
 &\ \ \ \ +\, \p_P g (u_p, \D u_p) : \D u_p \Big) \, \d {\mathcal{L}}^{n}\\
 &\leq \, C_8 \, {\av}_{\!\!\!\Om} g(u_p, \D u_p)^{p}\, \d \mL^n.
\end{split}
\]
The above two estimates, combined with the Euler-Lagrange equations, imply that $\lambda_p>0$. Hence, we may therefore define $\Lambda_p:=(\lambda_p)^{\frac{1}{p}}>0$. We will now obtain the upper and lower bounds. We determine the lower bound as follows:
\[
\begin{split}
 C_1(L_p)^p \, &=\, C_1 \, {\av}_{\!\!\!\Om} f(\D^2 u_p)^{p}\, \d \mL^n
 \\
 &\leq  \, {\av}_{\!\!\!\Om}f^{p-1}(\D^2u_p) \partial f(\D^2u_p): \D^2 u_p \, \d \mL^n
 \\
 &= \, \lambda_p \, {\av}_{\!\!\!\Om} g(u_p, \D u_p)^{p-1}\Big(\p_\eta g (u_p, \D u_p) \cdot \phi  \,  + \,  \p_P g (u_p, \D u_p) : \D u_p \Big) \, \, \d {\mathcal{L}}^{n}
 \\
 &\leq  \, \lambda_p C_8.
\end{split}
\]
Hence,
\[
\bigg(\frac{C_1}{C_8}\bigg)^{\frac{1}{p}}L_p\leq (\lambda_p)^{\frac{1}{p}}=\Lambda_p.
\]
The upper bound is determined analogously, by reversing the direction of the inequalities. Combining both bounds, we obtain the desired estimate.
\qed
\ms

\begin{proposition}\label{Proposition 6}
There exists $(u_{\infty}, \Lambda_{\infty})\in W^{2,\infty}_{\mathrm B}(\Om;\R^N)\times (0, \infty)$ such that, along a sequence $(p_j)_1^{\infty}$ of exponents, we have
\[
\left\{ \ \
\begin{array}{ll}
u_{p} \larrow u_\infty, &  \ \ \text{in } C^1 \big(\overline{\Om};\R^N \big), \smallskip
\\
\D^2 u_{p} \weak \D^2 u_\infty, &  \ \ \text{in } L^q (\Om; \R_s^{N\times n^2}), \ \text{for all} \ \ q\in (1,\infty), \smallskip
\\
\Lambda_p \larrow \Lambda_{\infty}, &  \ \ \text{in } [0,\infty),  \smallskip
\end{array}
\right.
\]
as $p_j\to\infty$. Additionally, $u_{\infty}$ solves the minimisation problem \eqref{1.1} and $\Lambda_{\infty}$ is given by \eqref{1.6}. Finally $\Lambda_{\infty}$ satisfies the uniform bounds \eqref{1.7}.
\end{proposition}

\BPP \ref{Proposition 6}.
Fix $p>n/\al, q\leq p$ and a map $v_0\in W^{2,\infty}_{\mathrm B}(\Om;\R^N)\setminus \{0\}$. Then, by Lemma \ref{lemma2} there exists $(t_p)_{p\in (n/\al, \infty]} \sub (0, \infty)$ such that $t_p \to t_{\infty}$ as $p\to \infty$ and satisfying $\|g(t_p v_0, t_p \D v_0)\|_{L^p(\Om)}=1$ for all $p\in (n/\al, \infty]$. By H\"older's inequality and minimality, we have the following estimate
\[
\begin{split}
\big\|f(\D^2u_p) \big\|_{L^q(\Om)} \, &\leq \,  \big\|f(\D^2u_p) \big\|_{L^p(\Om)}
\\
&\leq \,  \big\|f(t_p\D^2v_0) \big\|_{L^p(\Om)}
\\
& \leq\,   \big\| f(t_p\D^2v_0)  \big\|_{L^{\infty}(\Om)} 
\\
&\leq K\,+\,  \big\| f(t_\infty\D^2v_0)  \big\|_{L^{\infty}(\Om)}\\
& < \, \infty,
\end{split}
\]
for some $K>0$. By \eqref{1.3}(d), we have the bound $f^q(X)\geq C_4(q)|X|^{\al q}- C_3(q)$ for some constants $C_3(q), C_4(q)>0$ and all $X\in \R^{N \by n^2}_s.$ By the previous bound, we conclude that 
\[
\sup_{q\geq p}\|\D^2u_p\|_{L^{\al q}(\Om)}\leq C(q)< \infty,
\]
for some $q$-dependent constant. By arguing as in the proof of Lemma \ref{lemma3} through the use of Poincar\'e inequalities, we can conclude in both cases of boundary conditions with the bound
\[
\sup_{q\geq p}\|u_p\|_{W^{2, \al q}(\Om)}\leq C(q)< \infty,
\]
for a new $q$-dependent constant $C'(q)>0$. Standard compactness in Sobolev spaces and a diagonal sequence argument imply the existence of a mapping
\[ 
u_\infty \in \bigcap_{n/\al<p<\infty} W^{2, \al p}_{\mathrm B}(\Om; \R^N) 
\]
and a subsequence $(p_j)_1^\infty$ such that the desired modes of convergence hold true as $p_j\to \infty$ along this subsequence of indices. Fix a map $v\in W^{2, \infty}_{\mathrm B}(\Om; \R^N)$ satisfying the required constraint, namely $\|g(v, \D v) \|_{L^{\infty}(\Om)}=1$. In view of Lemma \ref{lemma2}, there exists $(t_p)_{p\in (n/\al, \infty)}\sub (0, \infty)$ satisfying that $t_p\to 1$ as $p\to \infty$, and additionally $\|g(t_p v, t_p \D v)\|_{L^p(\Om)}=1$ for all $p>n/\al$. By  H\"older's inequality, the definition of $L_p$ and minimality, we have 
\[
\big\|f(\D^2 u_p)\big\|_{L^q(\Om)} \, \leq \,  \big\|f(\D^2 u_p)\big\|_{L^p(\Om)} \, = \, L_p \, \leq  \, \big\|f(t_p\D^2 v)\big\|_{L^p(\Om)},
\]
for any such $v$. By the weak lower semi-continuity of the functional on $W^{2, \al q}_{\mathrm B}(\Om; \R^N)$, we may let $p_j \to\infty$ to obtain
\[
\begin{split}
\big\|f(\D^2u_{\infty})\big\|_{L^q(\Om)} \, &\leq \, \liminf_{p_j\to \infty}L_p
\\
&\leq \, \limsup_{p_j\to \infty} L_p 
\\
&\leq \, \limsup_{p_j \to \infty} \|f(t_{p_j}\D^2 v)\|_{L^p(\Om)}
\\
& =\, \|f(\D^2 v)\|_{L^{\infty}(\Om)}.
\end{split}
\] 
Now we may let $q\to \infty$ in the above bound, hence producing
\[
\big\|f(\D^2 u_{\infty})\big\|_{L^{\infty}(\Om)} \, \leq  \, \liminf_{p_j\to \infty}L_p  \, \leq \,  \limsup_{p_j\to \infty} L_p  \, \leq \, \|f(\D^2 v)\|_{L^{\infty}(\Om)}.
\]
for all mappings $v\in W^{2, \infty}_{\mathrm B}(\Om; \R^N)$ satisfying the constraint $\|g(v, \D v)\|_{L^{\infty}(\Om)}=1$. If we additionally show that in fact $u_{\infty}$ satisfies the constraint in \eqref{1.1}, then the above estimate shows both that it is the desired minimisers (by choosing $v:=u_\infty$), and also that the sequence $(L_{p_j})_1^\infty$ converges to the infimum. Now we show that this is indeed the case. In view of assumption \eqref{1.3}(d), the previous estimate implies also that $\D^2u_{\infty}\in \smash{L^{\infty}\big(\Om; \R_s^{N \by n^2}\big)}$, which together with Poincar\'e inequalities (as in the proof of Lemma \ref{lemma3}) shows that in fact $u_{\infty} \in W^{2, \infty}_{\mathrm B}(\Om; \R^N)$. By the continuity of the function $g$ and the fact that $u_{p} \larrow u_\infty$ in $C^1 \big(\overline{\Om};\R^N \big)$, we have
\[
\begin{split}
1&= \, \|g(u_p, \D u_p)\|_{L^p(\Om)}\\
&= \, \|g(u_{\infty}, \D u_{\infty})\|_{L^p(\Om)}\,+\, \|g(u_p, \D u_p)\|_{L^p(\Om)}-\|g(u_{\infty}, \D u_{\infty})\|_{L^p(\Om)}\\
&= \, \|g(u_{\infty}, \D u_{\infty})\|_{L^p(\Om)}\,+ \,\mathrm O\Big(\|g(u_p, \D u_p)- g(u_{\infty}, \D u_{\infty})\|_{L^{\infty}(\Om)}\Big)\\
&\!\!\larrow \|g(u_{\infty}, \D u_{\infty})\|_{L^{\infty}(\Om)},
\end{split}
\]
as $p_j\to \infty$. Consequently, $u_{\infty}$ satisfies the constraint, and therefore lies in the admissible class of \eqref{1.1}. Since $v$ was arbitrary in the energy bound, we conclude that $u_{\infty}$ in fact solves \eqref{1.1}. let us now define
\[
\Lambda_{\infty}:=\big\|f(\D^2 u_{\infty})\big\|_{L^{\infty}(\Om)}.
\]
We now show that $\Lambda_{\infty}>0$. Indeed, by our assumptions \eqref{1.3}-\eqref{1.4}, there is  only one map in $W^{2,\infty}(\Om;\R^N)$ satisfying $\|f(\D^2u_0)\|_{L^{\infty}(\Om)}=0$ and $|u_0|\equiv 0$ on $\p\Om$, namely the trivial map $u_0\equiv 0$, but $u_0$ is not contained in the admissible class of \eqref{1.1} because $\|g(u_0, \D u_0)\|_{L^{\infty}(\Om)}=0$. We now show that $\Lambda_p \larrow \Lambda_{\infty}$ as $p_j \to \infty$. By our earlier energy estimate, we have $L_p\larrow \Lambda_{\infty}$ as $p_j \to \infty$.
By Lemma \ref{lemma5}, we have
\[
0< \, \lim_{p_j\to \infty} \bigg(\frac{C_1}{C_8}\bigg)^{\frac{1}{p}}L_p  \, \leq  \,  \lim_{p_j\to \infty} \Lambda_p \, \leq  \, \lim_{p_j\to \infty} \bigg(\frac{C_2}{C_7}\bigg)^{\frac{1}{p}}L_p,
\]
and therefore $\Lambda_p \larrow \Lambda_{\infty}$ as $p_j \to \infty$.
To complete the proof we must deduce the claimed bounds for $\Lambda_{\infty}$. We first establish the lower bound. By utilising the Poincar\'e and Poincar\'e-Wirtinger inequalities (recall the proof of Lemma \ref{lemma3}) and that $g(0,0)=0$, we estimate
\[
\begin{split}
1&= \,   \big\|g(u_{\infty}, \D u_{\infty})  \big\|_{L^{\infty}(\Om)}
\\
&\leq \, \mathrm{diam}(\Om)  \big\|\D (g(u_{\infty}, \D u_{\infty})) \big\|_{L^{\infty}(\Om)}
\\
&\leq \, \mathrm{diam}(\Om)\Big(   \big\|\p_\eta g (u_{\infty}, \D u_{\infty})\D u_{\infty}  \big\|_{L^{\infty}(\Om)} \, +\,  \big\| \p_P g (u_{\infty}, \D u_{\infty})\D^2 u_{\infty}  \big\|_{L^{\infty}(\Om)} \Big)
\\
&\leq \, \mathrm{diam}(\Om)\Big(   \big\|\p_\eta g  \big\|_{L^{\infty} ((u_{\infty}, \D u_{\infty})(\overline{\Om}))}\|\D u_{\infty}\|_{L^{\infty}(\Om)} 
\\
&  \ \ \ \ \ \ \ \ \ \ \ \ \ \ \ \  + \|\p_P g \|_{L^{\infty} ((u_{\infty}, \D u_{\infty})(\overline{\Om} ))}\|\D^2 u_{\infty}\|_{L^{\infty}(\Om)}\Big)
\\
&\leq \, \| \D^2 u_{\infty}\|_{L^{\infty}(\Om)} \mathrm{diam}(\Om) \Big( C(\infty,\Om) \|\p_\eta g \|_{L^{\infty}((u_{\infty}, \D u_{\infty})(\overline{\Om} ))} 
\\
&  \ \ \ \ \ \ \ \ \ \ \ \ \ \ \ \ \ \ \ \ \ \ \ \ \ \ \ \ \ \ \ \ \ \ \ + \,  \|\p_P g \|_{L^{\infty}((u_{\infty}, \D u_{\infty})(\overline{\Om})}\Big),
\end{split}
\]
where $C(\infty,\Om)>0$ is the maximum of the Poincar\'e and the Poincar\'e-Wirtinger inequality constants on $\Om$ for $p=\infty$ (with the former being equal to $\mathrm{diam}(\Om)$). As $g\geq 0$ and $\|g(u_{\infty}, \D u_{\infty})\|_{L^{\infty}(\Om)}=1$, we conclude that $0\leq g(u_{\infty}, \D u_{\infty})\leq 1$ on $\overline{\Om}$. Hence
$(u_{\infty}, \D u_{\infty})(\overline{\Om}) \sub \{0\leq g \leq 1\}=\{g \leq 1\}$. Thus
\[
\begin{split}
1&\leq \| \D^2 u_{\infty}\|_{L^{\infty}(\Om)} \mathrm{diam}(\Om)  \Big( C(\infty,\Om) \|\p_\eta g \|_{L^{\infty}(\{ g \leq 1\})} \, +\, \|\p_P g \|_{L^{\infty}(\{g \leq 1\})}\Big)
\end{split}
\]
Rearranging assumption \eqref{1.3}(d), we may write $|X|\leq C_4^{-\frac{1}{\al}}(f(X)+C_3)^{\frac{1}{\al}}$, for any $X \in \R^{N \by n^2}_s$. Combining this inequality with the previous bound, we deduce
\[
\begin{split}
C_4^{\frac{1}{\al}}  \leq \,  \Big(\big\|f(\D^2 u_{\infty})\big\|_{L^{\infty}(\Om)} +\, C_3\Big)^{\!\frac{1}{\al}}  \mathrm{diam}(\Om)  \Big( & C(\infty,\Om) \|\p_\eta g \|_{L^{\infty}(\{ g \leq 1\})}  \\
& + \,  \|\p_P g \|_{L^{\infty}(\{g \leq 1\})}\Big),
\end{split}
\]
which leads directly to the claimed lower bound for the eigenvalue.

Now we establish the upper bound for $\La_\infty$. Since $\Om$ is by assumption a bounded domain with $C^2$ boundary, by standard results (see e.g.\ \cite[Sec.\ 14.6]{GT}), the distance function
\[
d_\Om \equiv \mathrm{dist}(\cdot,\p\Om) \ : \ \ \R^n \larrow \R,
\]
which is in $W^{1,\infty}_{\mathrm{loc}}(\R^n)$, is also $C^2$ on an inner tubular neighbourhood of $\p\Om$, namely there exists $\e_0 \in (0,1)$ such that
\[
d_\Om \in C^2(\overline{\Om^{\e_0}}), \ \ \ \Om^\e := \{d_\Om <\e\} \cap \Om.
\]
Let us also for convenience symbolise $\Om_\e := \{d_\Om >\e\} \cap \hspace{1pt} \Om$. Let us also fix $k\in\{1,2\}$, a unit vector $e\in \R^N$ and $\ze \in C^2(\R^n)$ with $\ze \equiv 0$ on $\Om_{\e_0}$. Then, for any $t_0>0$, the map $\xi_0 := t_0 (d_\Om)^k\ze e$ satisfies
\[
\xi_0 \, \in \, C^2(\overline{\Om};\R^N).
\]
Since $\d_\Om =0$ on $\p\Om$ and also $\D(\d^2_\Om)=0$ on $\p\Om$, it follows that $\xi_0 \in W^{2,\infty}_{\mathrm H}(\Om;\R^N)$ if $k=1$, whilst $\xi_0 \in W^{2,\infty}_{\mathrm C}(\Om;\R^N)$ if $k=2$. We will consider both cases simultaneously and declare this as 
\[
\xi_0 \in W^{2,\infty}_{\mathrm B}(\Om;\R^N).
\]
By Lemma \ref{lemma2}, we can adjust the constant $t_0>0$ to arrange
\[
\big\| g(\xi_0,\D \xi_0)\big\|_{L^\infty(\Om)} = 1.
\]
Hence, $\xi_0$ is in the admissible class of the minimisation problem \eqref{1.1}. By minimality and assumption \eqref{1.3}, we have the estimate
\beq
\label{2.3}
\La_\infty \, \leq\, C_5 (t_0)^\al \Big(\big\| \D^2(d_\Om^k \ze)\big\|_{L^\infty(\Om^{\e_0})}\Big)^\al  +\, C_6.
\eeq
By a direct computation, we have
\beq
\label{2.4}
\left\{ 
\begin{split}
\D^2(d_\Om^k \ze) \, &= \, k\big[ (k-1)\D d_\Om \ot \D d_\Om + d_\Om^{k-1}\D^2 d_\Om  \big] \ze \, +\, d_\Om^k \D^2 \ze
\\
&\ \  \ + \, k d_\Om^{k-1} \big( \D d_\Om\ot \D \ze  \, +\, \D \ze \ot \D d_\Om \big),
\end{split}
\right.
\eeq
on $\overline{\Om}$. For any $x\in \overline{\Om^{\e_0}}$, let us set $\mathrm{P}_{\Om}(x) := \mathrm{Proj}_{\p\Om}(x)$. Then, by \cite[Sec.\ 14.6, L.\ 14.17]{GT}, it follows that $|x- \mathrm{P}_{\Om}(x)|=d_\Om(x)$, and we also have the next estimates
\beq
\label{2.5}
\left\{  \ \ 
\begin{split}
\| d_\Om \|_{L^\infty(\Om^{\e_0})}  &\leq \e_0, \phantom{\Big|}
\\
 \| \D d_\Om \|_{L^\infty(\Om^{\e_0})} &\leq 1, 
\\
 \big\| \D^2 d_\Om \big\|_{L^\infty(\Om^{\e_0})} &\leq \sum_{i=1}^{n-1} \left\| \frac{ \ka_i \circ \mathrm{P}_{\Om} }{ 1 - ( \ka_i \circ \mathrm{P}_{\Om})d_\Om } \right\|_{L^\infty(\Om^{\e_0})},
\end{split}
\right.
\eeq
where $\{\ka_1,...,\ka_{n-1}\}$ are the principal curvatures of $\p\Om$. By \eqref{2.3}-\eqref{2.5} we have the estimate
\beq
\label{2.6}
\begin{split}
\La_\infty \, \leq\, C_5 (2t_0)^\al \Bigg( & \| \D^2 \ze \|_{L^\infty(\Om)} \, +\,  \| \D \ze \|_{L^\infty(\Om)}
\\
& +  \| \ze \|_{L^\infty(\Om)} \left( 1+ \sum_{i=1}^{n-1} \left\| \frac{ \ka_i \circ \mathrm{P}_{\Om} }{ 1 - ( \ka_i \circ \mathrm{P}_{\Om})d_\Om } \right\|_{L^\infty(\Om^{\e_0})} \right)  \!\!\Bigg)^{\!\al}  +\, C_6.
\end{split}
\eeq
It remains to select an appropriate function $\ze$ in order to estimate its derivatives in terms of the geometry of $\Om$, and to obtain an estimate for $t_0$. For the former, we argue as follows. Let $(\eta^\de)_{\de>0}$ be the family of standard mollifying kernels, as e.g.\ in \cite{KV}. We select
\[
\ze := \eta^{\e_0} * (\chi_{\R^n\set\Om}),
\]
which is the regularisation of the characteristic of the complement of $\Om$. It follows that this function satisfies the initial requirements, and additionally
\[
\left\{ 
\begin{split}
\D \ze &= \eta^{\e_0} * (\D \chi_{\R^n\set\Om}) = \eta^{\e_0} * \big(\mathcal{H}^{n-1}\LL_{\p\Om} \D \d_\Om \big),
\\
\D^2 \ze &= \D\eta^{\e_0} * (\D \chi_{\R^n\set\Om}) = \frac{1}{\e_0}(\D\eta)^{\e_0} * \big(\mathcal{H}^{n-1}\LL_{\p\Om} \D \d_\Om \big),
\end{split}
\right.
\]
by standard properties of convolutions and the differentiation of BV functions (see e.g.\ \cite{F} and \cite[Ch. 5, p. 198]{EG}). Then, by Young's inequality for convolutions, we have the estimates
\beq
\label{2.7}
\left\{ 
\begin{split}
 \| \D^2 \ze \|_{L^\infty(\R^n)} & \leq  \frac{C}{\e_0^{n+1}}\mathcal{H}^{n-1}(\p\Om),
\\
\| \D \ze \|_{L^\infty(\R^n)} &\leq  \mathcal{H}^{n-1}(\p\Om),
\\
\| \ze \|_{L^\infty(\R^n)} &\leq 1, \phantom{\Big|}
\end{split}
\right.
\eeq
for some universal constant $C>0$. Now we work towards an estimate for $t_0$ appearing in \eqref{2.3}. By assumption \eqref{1.4}, we have that the sublevel sets $\{g\leq t\}$ are compact in $\R^N \by \R^{N\by n}$ for any $t\geq0$. Let us define $R(t)$ as the smallest radius of the $N$-dimensional ball, for which $\{g\leq t\}$ is contained into the cylinder $\bar\mB^N_{R(t)}(0) \by \R^{N\by n}$:
\beq
\label{2.8}
R(t) \,:=\, \inf \Big\{ R>0 \, : \, \{g\leq t\} \sub \mB^N_R(0) \by \R^{N\by n}\Big\}.
\eeq
Then, we define a strictly increasing function $\rho : [0,\infty) \larrow [0,\infty)$ by setting
\beq
\label{2.9}
\rho(t) \,:=\, t+\sup_{0\leq s\leq t}R(s).
\eeq
Then, $\rho$ satisfies $\rho(0)=0$, and also that 
\[
\{g\leq t\} \, \sub \, \bar\mB^N_{\rho(t)}(0) \by \R^{N\by n},
\]
for any $t\geq0$. Further, by construction,
\[
\Big\{(\eta, P) \in \R^N \by \R^{N\by n}\, :\ \rho^{-1}(|\eta|)\leq t \Big\}\, =\, \bar\mB^N_{\rho(t)}(0) \by \R^{N\by n}.
\]
The above imply
\[
\rho^{-1}(|\eta|) \,\leq\, g(\eta,P),\ \ \ (\eta, P) \in \R^N \by \R^{N\by n}.
\]
Next, since $d_\Om^k\ze$ vanishes on $\p\Om \cup \overline{\Om_{\e_0}}$ and $g(0,0)=0$, we have
\[
\begin{split}
1\, &=\, \|g(\xi_0,\D\xi_0)\|_{L^\infty(\Om)}
\\
&=\, \sup_{\Om^{\e_0}}g(\xi_0,\D\xi_0)
\\
&\geq\, \sup_{\Om^{\e_0}}\rho^{-1}(|\xi_0|)
\\
&\geq\, \sup_{\Om^{\e_0}}\rho^{-1}\big(t_0 |d_\Om^k \ze|\big).
\end{split}
\]
Since $d_\Om\equiv \e_0/4$ on $\p\Om^{\e_0/4}$, and $\rho^{-1}$ is strictly increasing, the above implies
\beq
\label{2.10}
\begin{split}
1\, &\geq\, \max_{\p\Om^{\e_0/4}}\rho^{-1}\big(t_0 |d_\Om^k \ze|\big)
\\
&=\, \max_{\p\Om^{\e_0/4}}\rho^{-1}\Big(t_0 \Big(\frac{\e_0}{4}\Big)^k \ze\Big)
\\
&=\, \rho^{-1}\Big(t_0 \Big(\frac{\e_0}{4}\Big)^k \max_{\p\Om^{\e_0/4}}\ze\Big).
\end{split}
\eeq
Now we estimate $\max_{\p\Om^{\e_0/4}}\ze$ from below. Fix $x\in \p\Om^{\e_0/4}$. Then, since the standard mollifying kernel $\eta$ is a radial function (see e.g.\ \cite{KV}), there exists a universal $c>0$ such that $\eta\geq c$ on $\mB_{1/2}(0)$. Therefore,
\[
\begin{split}
\ze(x) \, &=\, \frac{1}{\e_0^n}\int_{\mB_{\e_0}(x)} \chi_{\R^n \set\Om}  \eta\Big( \frac{|y-x|}{\e_0}\Big) \mathrm d y
\\
&\geq\, \frac{1}{\e_0^n}\int_{\mB_{\e_0/2}(x) \set\Om}  \eta\Big( \frac{|y-x|}{\e_0}\Big) \mathrm d y
\\
&\geq\, \frac{c}{\e_0^n}\mL^n\big(\mB_{\e_0/2}(x) \set\Om \big),
\end{split}
\]
for any $x\in \p\Om^{\e_0/4}$. Finally, since $\p\Om$ satisfies the exterior sphere condition, the set $\mB_{\e_0/2}(x) \set\Om$ contains a ball $\mB_r(\bar x)$ centred at some point $\bar x$, where the maximum possible radius $\bar r$ is given by
\[
\bar r\, =\, \min\bigg\{\frac{\e_0}{8} \,,\, \underset{i=1,...,n-1}{\min}\frac{1}{\|\ka_i\|_{C^0(\p\Om)}}  \bigg\} .
\]
Therefore, if $\om(n)$ is the volume of the unit ball in $\R^n$,
\[
\begin{split}
\ze(x) \, &\geq\, \frac{c}{\e_0^n}\mL^n\big(\mB_{\e_0/2}(x) \set\Om \big) 
\\
&\geq\, \frac{c}{\e_0^n}\mL^n(\mB_{\bar r}(\bar x)) 
\\
& =\,  \frac{c}{\e_0^n}\om(n) \bar r^n
\\
& =\,  \frac{c \om(n)}{\e_0^n}  \min\bigg\{\Big(\frac{\e_0}{8}\Big)^n \,,\, \underset{i=1,...,n-1}{\min}\frac{1}{\big(\|\ka_i\|_{C^0(\p\Om)}\big)^n}  \bigg\}
\\
& =\,  c \om(n)   \min\bigg\{\frac{1}{2^{3n}} \,,\, \underset{i=1,...,n-1}{\min}\frac{1}{\big(\e_0\|\ka_i\|_{C^0(\p\Om)}\big)^n}  \bigg\},
\end{split}
\]
for any $x\in \p\Om^{\e_0/4}$. Hence, we have established the lower bound
\beq
\label{2.11}
\max_{\p\Om^{\e_0/4}}\ze \, \geq\,  c \om(n) \min\bigg\{\frac{1}{2^{3n}} \,,\, \underset{i=1,...,n-1}{\min}\frac{1}{\big(\e_0\|\ka_i\|_{C^0(\p\Om)}\big)^n}  \bigg\}.
\eeq
By \eqref{2.10} and \eqref{2.11}, we infer (since $\e_0<1$ and $k\in\{1,2\}$) that
\beq
\label{2.12}
\begin{split}
t_0 \, &\leq\,  \frac{4^k \rho(1)\e_0^{n-k}}{c \om(n)}\frac{1}{\min\bigg\{\dfrac{1}{2^{3n}}\, , \,  \underset{i=1,...,n-1}{\min}\dfrac{1}{\big(\e_0\|\ka_i\|_{C^0(\p\Om)}\big)^n}  \bigg\}}
\\
&\leq\,  \frac{32 \rho(1)}{c \om(n)}\Big(2^{3n} + \underset{i=1,...,n-1}{\max}\big(\|\ka_i\|_{C^0(\p\Om)}\big)^n  \Big).
\end{split}
\eeq
By \eqref{2.6}, \eqref{2.7}, and \eqref{2.12}, we conclude with the following upper bound for the eigenvalue:
\beq
\label{2.13}
\begin{split}
\La_\infty \, \leq &\ C_6 \, + \, C_5 \left[  \frac{16\rho(1)}{c \om(n)}\Big(2^{3n} \,+ \, \underset{i=1,...,n-1}{\max}\big(\|\ka_i\|_{C^0(\p\Om)}\big)^n  \Big) \right]^\al \centerdot
\\
&\centerdot \Bigg\{ 1\,+\, \bigg(1+\frac{C}{\e_0^{n+1}} \bigg) \mH^{n-1}(\p\Om) \  +  \, \sum_{i=1}^{n-1} \left\| \frac{ \ka_i \circ \mathrm{P}_{\Om} }{ 1 - ( \ka_i \circ \mathrm{P}_{\Om})d_\Om } \right\|_{L^\infty(\Om^{\e_0})}  \!\!\Bigg\}^{\!\al}.
\end{split}
\eeq
The claimed estimate \eqref{1.8A} follows from \eqref{2.13} above, by recalling that in view of \eqref{2.8}-\eqref{2.9}, we have
\[
\rho(1) \, =\, 1 \,+\, \sup_{0\leq t \leq 1}R(t),      
\]
and also that the last term of \eqref{2.13} is finite at least when
\[
\e_0 \, <\, \frac{1}{\underset{i=1,...,n-1}{\max} \|\ka_i\|_{C^0(\p\Om)}}.
\]
The result ensues. \qed
\ms

\begin{lemma}\label{lemma7}
For any $p>(n/\al) +2$, there exist measures $\nu_{\infty}\in \mathcal{M}(\overline{\Om})$ and ${\mathrm M}_\infty \in \mathcal{M}(\overline{\Om}; \R^{N \by n^2}_s)$ such that, along perhaps a further sequence $(p_j)_1^{\infty}$ of exponents, we have
\[
\left\{ \ \
\begin{array}{ll}
\nu_{p} \weakstar \nu_\infty, &  \ \ \text{in } \mathcal{M}(\overline{\Om}), \smallskip
\\
\mathrm M_{p} \weakstar {\mathrm M}_\infty, &  \ \ \text{in } \mathcal{M}(\overline{\Om}; \R^{N \by n^2}_s), \smallskip
\end{array}
\right.
\]
as $j\to\infty$, where the approximating measures  $\nu_p,\mathrm M_p$ are given by \eqref{1.11}.
\end{lemma}

\BPL \ref{lemma7}.
We begin by noting that since $g\geq 0$ and $\|g(u_p, \D u_p)\|_{L^p(\Om)}=1$, in view of \eqref{1.11} we have the bound
\[
\|\nu_p\|(\overline{\Om}) \, =\,\nu_p(\overline{\Om}) \, =\, {\av}_{\!\!\!\Om} g(u_p, \D u_p)^{p-1}\, \d \mL^n\leq \bigg({\av}_{\!\!\!\Om} g(u_p, \D u_p)^{p}\, \d \mL^n \bigg)^{\frac{p-1}{p}}=1.
\]
By the sequential weak$^*$ compactness of the space of Radon measures we can conclude that $\nu_{p} \weakstar \nu_\infty, \text{in } \mathcal{M}(\overline{\Om})$ up to the passage to a further subsequence. Now we establish appropriate total variation bounds for the measure $\mathrm M_p$. Since $f\geq0$, by the bounds of Lemma \ref{lemma5} and assumption \eqref{1.3}, we estimate (for sufficiently large $p$)
\[
\begin{split}
\|{\mathrm M}_p\|(\overline{\Om})\, &= \,  {\av}_{\!\!\!\Om} \bigg( \frac{f(\D^2 u_p)}{\Lambda_p}\bigg)^{p-1}|\partial f(\D^2 u_p)|  \, \d \mL^n
\\
&\leq  \, \frac{1}{\Lambda^{p-1}_p}  \, {\av}_{\!\!\!\Om} f(\D^2 u_p)^{p-1} \Big( C_5 f(\D^2 u_p)^{\beta}+C_6\Big) \, \d \mL^n
\\
&= \, \frac{C_5}{\Lambda^{p-1}_p}  \, {\av}_{\!\!\!\Om} f(\D^2 u_p)^{p-1+\beta}  \, \d \mL^n  \, + \,   \frac{C_6}{\Lambda^{p-1}_p}{\av}_{\!\!\!\Om} f(\D^2 u_p)^{p-1}  \, \d \mL^n.
\end{split}
\]
Hence,
\[
\begin{split}
\|{\mathrm M}_p\|(\overline{\Om})\, & \leq  \,  \frac{C_5}{\Lambda^{p-1}_p}\bigg({\av}_{\!\!\!\Om}f(\D^2 u_p)^p \, \d \mL^n \bigg)^{\frac{p-1+\beta}{p}}  \, + \, \frac{C_6}{\Lambda^{p-1}_p}\bigg({\av}_{\!\!\!\Om}f(\D^2 u_p)^{p}  \, \d \mL^n\bigg)^{\frac{p-1}{p}}
\\
&= \, C_5\frac{(L_p)^{p-1+\beta}}{\Lambda_p^{p-1}} \, + \, C_6\frac{(L_p)^{p-1}}{\Lambda_p^{p-1}}
\\
&= \, \bigg(\frac{L_p}{\Lambda_p}\bigg)^{p-1}\big(C_5L_p^{\beta}+C_6\big)
\\
&\leq  \, \bigg(\frac{C_8}{C_1}\bigg)^{1-\frac{1}{p}}\Big(C_5(\Lambda_{\infty}+1)^{\beta}+C_6 \Big).
\end{split}
\]
The above bound allows to conclude that $\mathrm M_{p} \weakstar {\mathrm M}_\infty \ \text{in} \ \mathcal{M}(\overline{\Om}; \R^{N \by n^2}_s)$, along perhaps a further subsequence of indices $(p_j)_1^\infty$ as $j\to\infty$. 
\qed
\ms

To conclude the proof of Theorem \ref{1}  we must ensure the PDE system \eqref{1.5} is indeed satisfied by the quadruple $(u_{\infty}, \Lambda_{\infty}, {\mathrm M}_\infty, \nu_{\infty})$.

\begin{lemma}\label{lemma8}
If ${\mathrm M}_\infty\in \mathcal{M}(\overline{\Om}; \R^{N \by n^2}_s)$ and $ \nu_{\infty}\in \mathcal{M}(\overline{\Om})$ are the measures obtained in Lemma \ref{lemma7}, then the pair $(u_{\infty}, \Lambda_{\infty})$ satisfies \eqref{1.5} for all $\phi\in C^2_{\mathrm B}(\overline{\Om};\R^N)$.
\end{lemma}

\BPL \ref{lemma8}.
Fix a test function $\phi \in C^2_{\mathrm B}(\overline{\Om};\R^N)$ and $p>n/\al+2$ by \eqref{1.11} we may rewrite the PDE system in \eqref{1.10} as follows
\[
 \int_{\Om} \D^2 \phi : \d \mathrm M_p
= \  \Lambda_p \, \int_{\Om} \Big( \p_\eta g (u_p, \D u_p) \cdot \phi + \p_P g (u_p, \D u_p) : \D \phi \Big)\, \d \nu_p,
\]
Recall that, by Proposition \ref{Proposition 6}, we have $\Lambda_p\larrow \Lambda_{\infty}$ and also $(u_p,\D u_p)\larrow (u_{\infty},\D u_{\infty})$ uniformly on $\overline{\Om}$, as $p_j\to \infty$. By assumption \eqref{1.4}(a), we have that $\p_\eta g (u_p,\D u_p)\larrow \p_\eta g (u_{\infty}, \D u_{\infty})$ and also $\p_P g (u_p,\D u_p)\larrow \p_P g (u_{\infty}, \D u_{\infty})$, both uniformly on $\overline{\Om}$, as $p_j\to \infty$. The result ensues by invoking Lemma \ref{lemma7}, in conjunction with weak$^*$-strong continuity of the duality pairing $\mathcal{M}(\overline{\Om})\times C(\overline{\Om}) \larrow \R$.

\subsection*{Acknowledgement} The authors would like to thank the referee of this paper for their careful reading of the manuscript, as well as their constructive suggestions.

\bibliographystyle{amsplain}

\end{document}